# The Approximate Bilinear Algorithm of Length 46 for Multiplication of 4 x 4 Matrices

**A. V. Smirnov**
S.Alexey.V@rambler.ru

**Abstract** — We propose the arbitrary precision approximate (APA) bilinear algorithm of length 46 for multiplication of 4 x 4 and 4 x 4 matrices. The algorithm has polynomial order 3 and 352 nonzero coefficients from total 2208.

## INTRODUCTION

In 1969, an unexpected result of Strassen directed attention of researchers to the problem of fast multiplication of large matrices (see [1]). Despite the efforts of numerous mathematicians, it is not yet solved in 2014. This problem in theory and practice requires constructing the bilinear algorithms with little length. The reader can find the introduction to the problem with comprehensive bibliography in the survey [2]. We give only one indispensable and not the most general definition. This concept of APA algorithms for matrix multiplication was introduced by Bini et al. [3]

Let $\gamma_{ij}^{t}$, $\alpha_{ij}^{t}$, and $\beta_{ij}^{t}$ be Laurent polynomials on some field, and let $a_{ij}$, $b_{ij}$ be elements of arbitrary ring.

**Definition**. A set of *3T* polynomial matrices

$$\gamma_{i_1 i_2}^{t}, \ i_1 = 1,2,...,n_1, \ i_2 = 1,2,...,n_3, \ t = 1,2,...,T,$$

$$\alpha_{j_1 j_2}^{t}, \ j_1 = 1,2,...,n_1, \ j_2 = 1,2,...,n_2, \ t = 1,2,...,T,$$

$$\beta_{k_1 k_2}^{t}, \ k_1 = 1,2,...,n_2, \ k_2 = 1,2,...,n_3, \ t = 1,2,...,T,$$

is an approximate (APA) bilinear algorithm of length $T$ for multiplication of matrices $n_1$ x $n_2$ and $n_2$ x $n_3$ if and only if

$$\sum_{t=1}^{T} \gamma_{i_1 i_2}^{t} \left( \sum_{j_1, j_2=1}^{n_1, n_2} \alpha_{j_1 j_2}^{t} a_{j_1 j_2} \right) \left( \sum_{k_1, k_2=1}^{n_2, n_3} \beta_{k_1 k_2}^{t} b_{k_1 k_2} \right) = \sum_{j=1}^{n_2} a_{i_1 j} b_{j i_2} + O(x).$$

The bilinear algorithm is noncommutative and can decrease the number of multiplication in the ring of matrix elements, needed for matrix multiplication. The $T$ can be less than $n_1 n_2 n_3$, if all the three dimensions are greater than one (see [1]).

The coefficients $\gamma_{ij}^{t}$, $\alpha_{ij}^{t}$, and $\beta_{ij}^{t}$ of approximate bilinear algorithm are the roots of Brent's system of cubic equations (see [4]).

$$\sum_{t=1}^{T} \gamma^{t}_{i_1 i_2} \alpha^{t}_{j_1 j_2} \beta^{t}_{k_1 k_2} = L_{i_1 i_2 j_1 j_2 k_1 k_2} + O(x),$$

where

$$L_{i_1 i_2 j_1 j_2 k_1 k_2} = \delta_{i_1 j_1} \delta_{j_2 k_1} \delta_{i_2 k_2},$$

$$\delta_{ij} = \begin{cases} 0, & i \neq j, \\ 1, & i = j. \end{cases}$$

Replacing each paired index in this system by a single one, we can write the objective function for this system in the form

$$S = \sum_{i=1}^{N_1} \sum_{j=1}^{N_2} \sum_{k=1}^{N_3} \left( \sum_{t=1}^{T} \gamma^{t}_{i} \alpha^{t}_{j} \beta^{t}_{k} - L_{ijk} \right)^2,$$

where $N_1 = n_1 n_3$, $N_2 = n_1 n_2$, $N_3 = n_2 n_3$.

This function allows to estimate the proper calculating error of the approximate algorithm. It also allows to use optimization methods for solution of Brent's system.

The approximate bilinear algorithm for multiplication of matrices 4 x 4 of length 47 was proposed in [5]. The existence of algorithm with length 46 and the partly heuristic computational optimization method for deriving bilinear algorithms of matrix multiplication were declared in [6].

Since all the coefficients of polynomials in the proposed algorithm are in {-1, 1}, we present the algorithm as the values of polynomials of argument $x = 0.1$. This form is excessive in data, but more convenient for computer output and input. We assume that research and using of this algorithm without computer are not easy. Beware that as above we use the traditional order of indexation in result matrices $\gamma^{t}$.

## ALGORITHM

| $\gamma^t$ | | | | $\alpha^t$ | | | | $\beta^t$ | | | |
|---|---|---|---|---|---|---|---|---|---|---|---|

$t = 1$

| | | | | | | | | | | | |
|---|---|---|---|---|---|---|---|---|---|---|---|
| 0 | 0 | 1 | 0 | 0.1 | 0.1 | −0.1 | 10 | 0 | 0 | 10 | 0 |
| 0.1 | 0 | 0 | 0 | 0 | 10 | 0 | 0 | 1 | 0 | 0 | 0 |
| 0 | 0 | 0 | 0 | 0 | 0 | 0 | 0 | 1 | 0 | 0 | 0 |
| 0 | 0 | 0 | 0 | 0 | 0 | 0 | 0 | 0 | 0 | 0.1 | 0 |

$t = 2$

| | | | | | | | | | | | |
|---|---|---|---|---|---|---|---|---|---|---|---|
| 0 | 0 | 0 | 0 | 0 | 0 | 0 | 0 | 0 | −10 | 0 | 10 |
| 0 | 0.1 | 0 | 0 | 0 | 0 | 1 | 10 | 0 | 0 | 0 | 0 |
| 0 | 0 | 0 | 10 | 0 | 0 | −0.1 | 0 | 0 | 0 | 0 | 0 |
| 0 | −0.1 | 0 | 0 | 0 | 0 | 0 | 0 | 0 | 1 | 0 | 0 |

$t = 3$

| | | | | | | | | | | | |
|---|---|---|---|---|---|---|---|---|---|---|---|
| 0 | 0 | 0 | 0 | 0 | 0 | 0 | 0 | 0 | 10 | 0.1 | 0 |
| 0 | 0 | 0 | 0 | 0 | 0 | 0 | 0 | 0 | 0 | 0 | 0 |
| 0 | −10 | 0 | 0 | −0.01 | 0 | 0 | −0.1 | 0 | 0 | 0 | 0 |
| 0 | 0.1 | 10 | 0 | 0 | 0 | 0 | 10 | 0 | 0 | 0.01 | 0 |

$t = 4$

| | | | | | | | | | | | |
|---|---|---|---|---|---|---|---|---|---|---|---|
| 0 | 0 | 0 | 0 | 0 | 0 | 0 | 0 | 0 | −10 | 0 | 0 |
| 0 | 0 | 0 | 0 | 0 | 0 | 1 | 10 | 0 | 0 | 0 | 0 |
| 0 | 0 | 0 | −10 | 0 | 0 | −0.1 | 10 | 0 | 0 | 0 | 0 |
| 0 | 0.1 | 0 | 0 | −1 | 0 | 0 | 10 | 0 | 1 | 0 | 0 |

$t = 5$

| | | | | | | | | | | | |
|---|---|---|---|---|---|---|---|---|---|---|---|
| 0 | 0 | 0 | 0 | 0 | 0 | 0 | 0 | 0 | 0 | 0 | 10 |
| 0 | 0 | 0 | 0 | 0 | 0 | 0 | 0 | 0 | 0 | 0 | 0 |
| 0 | 0 | 0 | 10 | 0.01 | 0 | 0.1 | 0 | 0 | 0 | 0 | 1 |
| 0 | 0 | 0 | −0.1 | −1 | 0 | 0 | 10 | 0 | 0 | 0 | 0 |

$t = 6$

| | | | | | | | | | | | |
|---|---|---|---|---|---|---|---|---|---|---|---|
| 0 | 0 | 0 | 0 | 0 | 0 | 0 | 0 | 0 | 0 | 0 | 0 |
| 0 | 0 | 0.1 | 10 | 0.01 | 10 | 1 | 0 | 0 | 0 | 0 | 0.01 |
| 0 | 0 | 0 | 0 | 0 | 0 | 0.01 | 0 | 0 | 0 | 10 | 0 |
| 0 | 0 | 0 | 0 | 0 | 0 | 0 | 0 | 0 | 0 | 0 | 0 |

$t = 7$

| | | | | | | | | | | | |
|---|---|---|---|---|---|---|---|---|---|---|---|
| 0 | 0 | 0 | 0 | 0 | 0 | 0 | 0 | 0 | 10 | 0 | 0 |
| 0 | 10 | 0 | 0 | 0.01 | 0 | 0.01 | 0.1 | 0 | 0 | 0.01 | 0 |
| 0 | 0 | 0 | 0 | 0 | 0 | 0 | 0 | 0 | 0 | 0 | 0 |
| 0 | 0 | 10 | 0 | 0 | 10 | 0 | 0 | 0 | 0 | 0 | 0 |

$t = 8$

| | | | | | | | | | | | |
|---|---|---|---|---|---|---|---|---|---|---|---|
| 0 | 0 | 0.1 | 10 | 0 | 0 | 0 | 10 | 0 | 0 | 0 | 0 |
| 0.01 | 0 | 0 | 0 | 0 | 0 | 0 | 0 | 0 | 0 | 0 | 0 |
| −10 | 0 | 0 | 0 | 0 | 0 | 0.01 | 0 | −10 | 0 | 0 | 0 |
| 0 | 0 | 0 | 0 | 0 | 0 | 0 | 0 | 0 | 0 | 0 | 0.01 |

$t = 9$

| | | | | | | | | | | | |
|---|---|---|---|---|---|---|---|---|---|---|---|
| 0 | 0.1 | 0 | 0 | 0 | 0 | 0 | 0 | 0 | 10 | 0 | 0 |
| 0 | 0 | 0 | 0 | 0 | 0 | 0 | 0 | 0 | 0 | 0 | 0 |
| 1 | −0.1 | 0 | −10 | 0 | 0 | 0 | 10 | 0 | 0 | 0 | 0 |
| 0 | 0.1 | 0 | 0 | 0 | 0 | 0 | 0 | 0 | −1 | 0 | −0.01 |

$t = 10$

| | | | | | | | | | | | |
|---|---|---|---|---|---|---|---|---|---|---|---|
| 0 | 0 | 0 | 0 | 0 | 0 | 0.1 | 0.1 | 0 | 0 | 0 | 0 |
| 0 | 0 | 0 | 0 | 0 | 0 | −0.1 | 0 | 0 | 0 | 0 | 0 |
| 0 | 0 | 0 | 0 | 0 | 0 | 0.1 | 0 | 0 | −1 | 0 | 0.01 |
| 0 | 0 | 0 | 10 | 0 | 10 | 10 | 10 | 0 | 0 | 0 | 0 |

$t = 11$

| | | | | | | | | | | | |
|---|---|---|---|---|---|---|---|---|---|---|---|
| 0 | 0 | 1 | 0 | 0 | 0 | 0 | −10 | 0 | 0 | 10 | 0 |
| 0.1 | 0 | 0 | 0 | 0 | 0 | 0 | 0 | 1 | 0 | 0 | 0 |
| 10 | 0 | 0 | 0 | 0 | 0.1 | 0 | 0 | 0 | 0 | 0 | 0 |
| 0 | 0 | 0 | 0 | 0 | 0 | 0 | 0 | 0 | 0 | 0 | 0 |

$t = 12$

| | | | | | | | | | | | |
|---|---|---|---|---|---|---|---|---|---|---|---|
| 0 | 0 | 0 | 0 | 0 | 0 | 0 | 0 | 0 | 0 | 0 | 0.1 |
| 0 | 0 | 0 | 0 | 0 | 0 | 0 | 0 | 0 | 0 | 0 | 0 |
| 0 | 10 | 0 | 0 | 0 | 0 | 0.1 | 0 | 0 | 1 | 0 | 0 |
| 0 | 0 | 0 | 10 | 0 | 0 | 0 | 10 | 0 | 0 | 0 | 0.01 |

$t = 13$

| | | | | | | | | | | | |
|---|---|---|---|---|---|---|---|---|---|---|---|
| 10 | 0 | 0 | 0 | 0 | 10 | 0 | 0 | -10 | 0 | 0 | 0 |
| 0 | 0 | 0 | 0 | 0 | 0 | 0 | 0 | 0.01 | 0 | 0 | 0 |
| 0 | 0 | 0 | 0 | 0 | 0 | 0 | 0 | 0 | 0 | 0 | 0 |
| -1 | 0 | 0 | 0 | 0.1 | 0 | 0 | 0 | 0.01 | 0 | 0 | 0 |

$t = 14$

| | | | | | | | | | | | |
|---|---|---|---|---|---|---|---|---|---|---|---|
| 0 | 0 | 0 | 10 | 0 | 0 | 0 | 0 | 0 | 0 | 0 | -10 |
| -0.1 | 0.1 | 0 | 10 | 0 | 0 | 1 | 10 | 0 | 0 | 0 | 0 |
| 0 | 0 | 0 | 10 | 0 | 0 | 0 | 0 | 0 | 0 | 0 | 0 |
| 0 | -0.1 | 0 | 0 | 0 | 0 | 0 | 0 | 0 | 0 | 0 | 0 |

$t = 15$

| | | | | | | | | | | | |
|---|---|---|---|---|---|---|---|---|---|---|---|
| 0 | 0 | 0 | 0 | 0 | 0 | 0 | 0 | 0 | 0 | 0 | 10 |
| 0 | 0 | 0 | 10 | 0.01 | 0 | 1 | 10 | 0 | 0 | 0 | 0 |
| 0 | 0 | 0 | 0 | 0 | 0 | 0 | 0 | 0 | 0 | -10 | 0.1 |
| 0 | 0 | 0 | 0 | 0 | 0 | 0 | 0 | 0 | 0 | 0 | 0.01 |

$t = 16$

| | | | | | | | | | | | |
|---|---|---|---|---|---|---|---|---|---|---|---|
| 0 | 1 | 0 | 0 | 0 | 0 | 0 | 10 | 0 | -1 | 0 | 0 |
| 0 | 0 | 0 | 0 | 0 | 0 | 0 | 0 | 0 | 0 | 0 | 0 |
| 10 | 0 | 0 | 0 | 0 | 0 | 0 | 10 | 0 | 0 | 0 | 0 |
| 0 | 0 | 0 | 0 | 0 | 0 | 0 | 0 | 0.01 | 0.1 | 0 | 0 |

$t = 17$

| | | | | | | | | | | | |
|---|---|---|---|---|---|---|---|---|---|---|---|
| 0 | 0 | -1 | 0 | 0 | 0 | -0.1 | 0 | 0 | 0 | 0 | 0 |
| 0 | 0 | 0.01 | 0 | 0 | 10 | 0 | 0 | 0 | 0 | 10 | 0 |
| 0 | 0 | 0 | 0 | 0 | 0 | 0 | 0 | 1 | 0 | 10 | 0 |
| 0 | 0 | 0 | 0 | 0 | 0 | 0 | 0 | 0 | 0 | 0 | 0 |

$t = 18$

| | | | | | | | | | | | |
|---|---|---|---|---|---|---|---|---|---|---|---|
| 0 | 0 | 0 | 0 | 0 | 0 | 0 | 0 | 0 | 0 | 0 | 0 |
| 0 | 10 | 0 | 0 | 0 | 0 | 0.1 | 0 | 0 | 0 | 0 | 0.01 |
| 0 | 0 | 0 | 0 | 0 | 0 | 0 | 0 | 0 | 1 | 0 | -0.01 |
| 0 | 0 | 0 | -10 | 0 | -10 | 0 | 0 | 0 | 0 | 0 | 0 |

$t = 19$

| | | | | | | | | | | | |
|---|---|---|---|---|---|---|---|---|---|---|---|
| 0 | 0 | 0 | 0 | 0 | 0 | 0 | 0 | 0 | 0 | 10 | 0 |
| -1 | 0 | 0 | 0 | 0 | 0 | 0 | -10 | 0 | 0 | 0.01 | 0 |
| 0 | 0 | 10 | 0 | 0.01 | 10 | 0 | 0 | 0 | 0 | 0 | 0 |
| 0 | 0 | 0 | 0 | 0 | 0 | 0 | 0 | 0 | 0 | 0 | 0 |

$t = 20$

| | | | | | | | | | | | |
|---|---|---|---|---|---|---|---|---|---|---|---|
| 0 | 0 | 0 | -10 | 0.01 | 0 | 0 | 0 | 0 | 0 | 0 | -10 |
| 0.1 | 0 | 0 | 0 | 0 | 0 | 1 | 10 | 0 | 0 | 0 | 0 |
| 0 | 0 | 0 | 0 | 0 | 0 | 0 | 0 | 10 | 0 | 0 | 0 |
| 0 | 0 | 0 | 0 | 0 | 0 | 0 | 0 | 0 | 0 | 0 | 0 |

$t = 21$

| | | | | | | | | | | | |
|---|---|---|---|---|---|---|---|---|---|---|---|
| 0 | 0 | 0 | 0 | 0 | 0 | 0 | 0 | 0 | 0 | 10 | 0 |
| 0 | 0 | 0 | 0 | 0 | 0 | 0 | 0 | 0 | 0 | 0 | -0.01 |
| 0 | 0 | 0 | -10 | 0 | 10 | 0 | 0 | 0 | 0 | 0 | 0 |
| 0 | 0 | 0.1 | 0 | 1 | 0 | 0 | -10 | 0 | 0 | 0 | 0 |

$t = 22$

| | | | | | | | | | | | |
|---|---|---|---|---|---|---|---|---|---|---|---|
| 0 | 0 | 0 | 0 | 0 | 0 | 0 | 0 | 0 | 0 | 0 | 0 |
| 0 | -10 | 0 | 0 | 0 | -0.01 | 0 | 0 | 0.01 | 10 | 0 | 0 |
| 0 | 0 | 0 | 0 | 0 | 0 | 0 | 0 | 0 | 0 | 0 | 0 |
| 10 | 0.01 | 0 | 0 | 0 | 10 | 0 | 0 | 0 | 0 | 0 | 0 |

$t = 23$

| | | | | | | | | | | | |
|---|---|---|---|---|---|---|---|---|---|---|---|
| 0 | 0 | 0 | 0 | 0 | 0 | 0 | 0 | 0 | 0 | 10 | 0 |
| 1 | 0 | 0 | 0 | 10 | −1 | 0 | −10 | 0 | 0 | 0 | 0 |
| 0 | 0 | 0 | 0 | 0.01 | 0 | 0 | 0 | 1 | 0 | 0 | 0 |
| 0 | 0 | 0 | 0 | 0 | 0 | 0 | 0 | −0.1 | 0 | 0 | 0 |

$t = 24$

| | | | | | | | | | | | |
|---|---|---|---|---|---|---|---|---|---|---|---|
| 0 | 0 | 1 | 0 | 0 | −0.1 | 0.1 | 0 | 0 | 0 | 10 | 0 |
| 0 | 0 | 0 | 0 | 0 | −10 | 0 | 0 | 1 | 0 | −10 | 0 |
| 0 | 0 | 0 | 0 | 0 | 0 | 0 | 0 | 0 | 0 | 0 | 0 |
| 0 | 0 | 0 | 0 | 0 | 0 | 0 | 0 | 0 | 0 | 0.1 | 0 |

$t = 25$

| | | | | | | | | | | | |
|---|---|---|---|---|---|---|---|---|---|---|---|
| 0 | 0 | 0 | 0 | 0 | 0 | 0 | 0 | 0 | 10 | 10 | 10 |
| 0 | 0 | 0 | 0 | 0 | 0 | 0 | 0 | 0 | 0 | 0 | −0.01 |
| 0 | 0 | 0 | 10 | 0 | 0 | 0 | 0 | 0 | 0 | 0 | 1 |
| 0 | 0 | 0 | 0 | 1 | 0 | 0 | −10 | 0 | −1 | 0 | 0 |

$t = 26$

| | | | | | | | | | | | |
|---|---|---|---|---|---|---|---|---|---|---|---|
| 10 | 0 | 0 | 0 | 0.01 | 10 | 0 | −10 | 10 | 0 | 0 | 0 |
| 0 | 0 | 0 | 0 | 0 | 0 | 0 | 0 | 0 | 0 | 0 | 0 |
| 0 | 0 | 0 | 0 | −0.01 | 0 | 0 | 0 | 0 | 0 | 0 | 0 |
| 0 | 0 | 0 | 0 | 0.1 | 0 | 0 | 0 | −0.01 | 0 | 0 | 0 |

$t = 27$

| | | | | | | | | | | | |
|---|---|---|---|---|---|---|---|---|---|---|---|
| 0 | 10 | 0 | 0 | 0.01 | 0 | 0 | 0.1 | 0 | 10 | 0 | 0 |
| 0 | 0 | 0 | 0 | 0 | 0 | 0 | 0 | 0 | 0 | 0 | 0 |
| 0 | 0 | 0 | 0 | 0 | 0 | 0 | 0 | 0 | −1 | 0.01 | 0 |
| 0 | 0 | 0 | 0 | 0 | 0 | 10 | 0 | 0 | 0 | 0 | 0 |

$t = 28$

| | | | | | | | | | | | |
|---|---|---|---|---|---|---|---|---|---|---|---|
| 0 | 10 | 0 | 0 | 0 | 0 | 0.1 | 0.1 | 0 | 0 | 0 | 0 |
| 0 | 0 | 0 | 0 | 0 | 0 | 0 | 0 | 0 | 0 | 0 | 0 |
| 0 | 0 | 0 | 0 | 0 | 0 | 0 | 0 | 0 | 1 | 0 | 0 |
| 0 | 0.1 | 0 | 10 | 0 | 0 | 10 | 0 | 0 | 0 | 0 | 0 |

$t = 29$

| | | | | | | | | | | | |
|---|---|---|---|---|---|---|---|---|---|---|---|
| 0 | 0 | 0 | 0 | 0 | 0 | 0 | 0 | 0 | 0 | 0 | 0 |
| 0 | 0 | 0 | 0 | 0 | 0 | 0 | 0 | 0 | 10 | 0 | 0 |
| 0 | 10 | 0 | 0 | 0 | 0.01 | 0 | 0 | 0 | 0 | 0 | 0 |
| 10 | 0 | 0 | 0 | 0 | 0 | 0 | 10 | 0.01 | 0 | 0 | 0 |

$t = 30$

| | | | | | | | | | | | |
|---|---|---|---|---|---|---|---|---|---|---|---|
| 0 | 0 | 0 | 1 | 0 | −10 | 0 | 0 | 0 | 0 | 0 | 0 |
| 0 | 0 | 0 | 0 | 0 | 0 | 0 | 0 | 0 | 0 | 0 | −0.1 |
| 0 | 0 | 0 | 0 | 0 | 0 | 0 | 0 | 10 | 0 | 0 | 0 |
| 1 | 0 | 0 | 0 | 0 | 0 | 0.1 | 0 | 0 | 0 | 0 | 0 |

$t = 31$

| | | | | | | | | | | | |
|---|---|---|---|---|---|---|---|---|---|---|---|
| 0 | 0.01 | 0 | 0 | 0 | 10 | 0 | 0 | −1 | 0 | 0 | 0 |
| 0 | 0 | 0 | 0 | 0 | 0 | 0 | 0 | 0 | 10 | 0 | −0.01 |
| 0 | 0 | 0 | 0 | 0 | 0 | 0 | 0 | 1 | 0 | 0 | 0 |
| 10 | 0 | 0 | 0 | 0 | 0 | 0 | 0 | 0 | 0 | 0 | 0 |

$t = 32$

| | | | | | | | | | | | |
|---|---|---|---|---|---|---|---|---|---|---|---|
| 0 | 0 | 0 | 0 | 0 | 0 | 0 | 10 | −10 | 1 | 10 | 0 |
| 0 | 0 | 0 | 0 | 0 | 0 | 0 | 0 | 1 | 0 | 0 | 0 |
| 10 | 0 | 0 | 0 | 0 | 0 | 0 | 0 | −10 | 0 | 0 | 0 |
| 0 | 0 | 0 | 0 | 0 | 0 | 0 | 0 | −0.01 | −0.1 | 0 | 0.01 |

$t = 33$

| | | | | | | | | | | | |
|---|---|---|---|---|---|---|---|---|---|---|---|
| 0 | 0 | 0 | 0 | 0 | 0 | 0 | 0 | 0.1 | 0 | −10 | 0 |
| 1 | 0 | −0.01 | 0 | 10 | 0 | 0 | 0 | 0 | 0 | 0 | 0 |

| 0 | 0 | 0 | 0 | 0 | 0 | 0 | 0 | -1 | 0 | 0 | 0 |
| 0 | 0 | 0 | 0 | 0 | 0 | 0 | 0 | 0.1 | 0 | 0 | 0 |

$t = 34$

| 0 | 0 | 0 | 10 | 0.01 | 1 | -10 | 10 | 0 | 0 | 0 | 0 |
| 0 | 0 | 0 | 0 | 0 | 0 | 1 | 10 | 0 | 0 | 0 | 0 |
| 0 | 0 | 0 | 0 | 0 | 0 | 0.01 | 0 | 10 | 0 | 0 | 0 |
| 0 | 0 | 0 | 0 | 0 | 0 | -0.01 | 0 | 0 | 0 | 0 | 0 |

$t = 35$

| 0 | 0 | 0 | 0 | 0 | 0 | 0 | 0 | 0 | -10 | -0.1 | 0.1 |
| 0 | 0 | 0 | 0 | 0 | 0 | 0 | 0 | 0 | 10 | 0 | 0 |
| 0 | 10 | 0 | 0 | 0 | 0 | 0 | 0 | 0 | 1 | 0 | 0 |
| 0 | 0 | 0 | 0 | 0 | 0 | 0 | -10 | 0.01 | 0 | -0.01 | 0.01 |

$t = 36$

| 0 | 0 | 0 | 0 | 0 | 0 | 0 | 0 | 0 | 0 | 0 | 0 |
| 0 | 0 | 0 | 10 | 0 | 0 | 0 | 10 | 0 | 0 | 0 | 0 |
| 0 | 0 | -10 | 0 | 0 | 0 | -0.01 | 0 | 0 | 0 | 10 | -0.1 |
| 0 | 0 | 0 | 0 | 0 | 0 | 0 | 0 | 0 | 0 | 0 | 0 |

$t = 37$

| 0 | 0 | 0 | 0 | 0 | 0 | 0 | 0 | 0 | -10 | 0 | 0 |
| 0 | 0 | 0 | 0 | 0.01 | 0 | 0.01 | 0.1 | 0 | 0 | 0 | 0 |
| 0 | 0 | 0 | 0 | -0.01 | 0 | 0 | -0.1 | 0 | 0 | 0 | 0 |
| 0 | 0 | 10 | 0 | 0 | 10 | 10 | 10 | 0 | 0 | 0 | 0 |

$t = 38$

| 0 | 0 | 0 | 0 | 0 | 0 | 0 | 0 | 0 | -10 | 0 | 0 |
| 0 | 10 | 0 | 0 | 0 | 0 | 0 | 0 | 0.01 | 10 | -0.01 | 0.01 |
| 0 | 0 | 0 | 0 | 0 | 0 | 0 | 0 | 0 | 1 | 0 | -0.01 |
| 0 | 0 | 0 | 0 | 0 | 10 | 0 | 0 | 0 | 0 | 0 | 0 |

$t = 39$

| 0 | 0 | 0.01 | 0 | 0 | 0 | 0 | 0 | 0 | 0 | 10 | 0 |
| -1 | 0 | 0 | 0 | 0 | 0 | 0 | 0 | 0 | 0 | 0 | 0 |
| 0.1 | 0 | 10 | -10 | 0 | -10 | 0 | 0 | 0 | 0 | 0 | 0 |
| 0 | 0 | 0.1 | 0 | 0 | 0 | 0 | 0 | 0 | 0 | 0 | 0 |

$t = 40$

| 0 | 0 | 0 | 0 | 0 | 10 | 0 | 0 | 0 | 0 | 0 | 0 |
| 0 | 0 | 0 | 0 | 0 | -0.01 | 0 | 0 | 0 | -10 | 0 | 0 |
| 0 | 0 | 0 | 0 | 0 | 0.01 | 0 | 0 | 0 | 0 | 0 | 0 |
| 10 | 0 | 0 | 0 | 0 | 10 | 0 | 10 | 0 | 0 | 0 | 0 |

$t = 41$

| 0 | 0 | 0 | 0 | 0 | 0 | 0 | 0 | 0 | 0 | 10 | 0 |
| 0 | 0 | 0 | 0 | 0 | 0 | 0 | 10 | 0 | 0 | 0.01 | 0 |
| 0 | 0 | 10 | 0 | 0 | 0 | 0 | 0 | 0 | 0 | 10 | -0.1 |
| 0 | 0 | 0 | 0 | 0 | 0 | 0 | 0 | 0 | 0 | -10 | 0 |

$t = 42$

| 0 | 0 | 0 | 0 | 0 | 0 | 0 | 0 | 0 | 0 | 0 | 0 |
| 0 | 0 | 0.01 | 0 | 0 | 0 | 0 | 10 | 0 | 0 | 0 | 0 |
| 0 | 0 | 10 | 0 | 0 | 0 | 0 | 0.01 | 0 | 0 | 0 | 0 |
| 0 | 0 | 0 | 0 | 0 | 0 | 0 | 0 | 0 | 0 | 10 | 0 |

$t = 43$

| 10 | 0 | 0 | 0 | 0 | 0 | 0 | 10 | 10 | 0 | 0 | 0 |
| 0 | 0 | 0 | 0 | 0 | 0 | 0 | 0 | 0 | 0 | 0 | 0 |
| 10 | 0 | 0 | 0 | 0.01 | 0 | 0 | 0 | 0 | 0 | 0 | 0 |

| 0 | 0 | 0 | 0 | 0 | 0 | 0 | 0 | 0 | 0 | 0 | 0 |
|---|---|---|---|---|---|---|---|---|---|---|---|

$t = 44$

| 0 | −10 | 0 | 0 | 0 | 0 | 0 | 0 | 0 | 10 | 0 | 0 |
|---|---|---|---|---|---|---|---|---|---|---|---|
| 0 | 0 | 0 | 0 | 0 | 0 | 0 | 0 | 0 | 0 | 0 | 0 |
| 0 | 0 | 0 | 0 | 0 | 0 | 0 | 0 | 0 | 0 | 0.01 | 0 |
| 0 | 0 | 10 | 0 | 0 | 0 | 10 | 0 | 0 | 0 | 0 | 0 |

$t = 45$

| 0.01 | 0 | 0 | 10 | 0 | 0 | 10 | 0 | 0 | 0 | 0 | 0 |
|---|---|---|---|---|---|---|---|---|---|---|---|
| 0 | 0 | 0 | 0 | 0 | 0 | 0 | 0 | 0 | 0 | 0 | 0 |
| 0 | 0 | 0 | 0 | 0 | 0 | 0 | 0 | 10 | 0 | 0 | 0.01 |
| 0 | 0 | 0 | 0 | 0 | 0 | 0 | 0 | 0 | 0 | 0 | 0 |

$t = 46$

| 0 | 0 | −1 | 0 | 0 | 0 | 0 | 0 | 0 | 0 | 0 | 0 |
|---|---|---|---|---|---|---|---|---|---|---|---|
| 0 | 0 | 0.11 | 10 | 0 | −10 | 0 | 0 | 0 | 0 | 0 | 0 |
| 0 | 0 | 0 | 0 | 0 | 0 | 0 | 0 | 0 | 0 | 10 | 0 |
| 0 | 0 | 0 | 0 | 0 | 0 | 0 | 0 | 0 | 0 | 0 | 0 |

## DISCUSSION

This algorithm is of polynomial order three, while the simplest approximate algorithms have polynomial order one (see [3]). The algorithm has 352 nonzero coefficients from total 2208. All coefficients except binomial $\gamma^{46}_{2,3}$ are monomials. All coefficients of polynomials are in {-1, 1}. The number of additions can be reduced by trivial eliminating duplicate calculations in linear forms. Objective function of the algorithm is

$$S = 37\,x^2 + 53\,x^4 + 37\,x^6 + O(\,x^8\,).$$

Our goal was to simplify the algorithm, therefore we deleted from it all monomials $\pm x^3$. That increased the first term of $S$.

We have the set of numerical coefficients in [-10, 10] with $S < 0.002$, but the bound of precision for approximate algorithm mostly depends on polynomial order of the algorithm. For calculations with 16 decimal digits of precision, the proposed polynomial algorithm achieves minimal total error ($S(x)$ and round-off error) for $x$ near 2E-5 and $S(x)$ near 1.48E-8. In this case, one recursion of the algorithm gives 4 precise digits in the result.